\documentstyle[amstex,epsf]{article}

\begin{document}
\def\qed{$\square$}

\def\crsald{\vbox to 15pt{\epsfxsize=25pt
\epsfbox{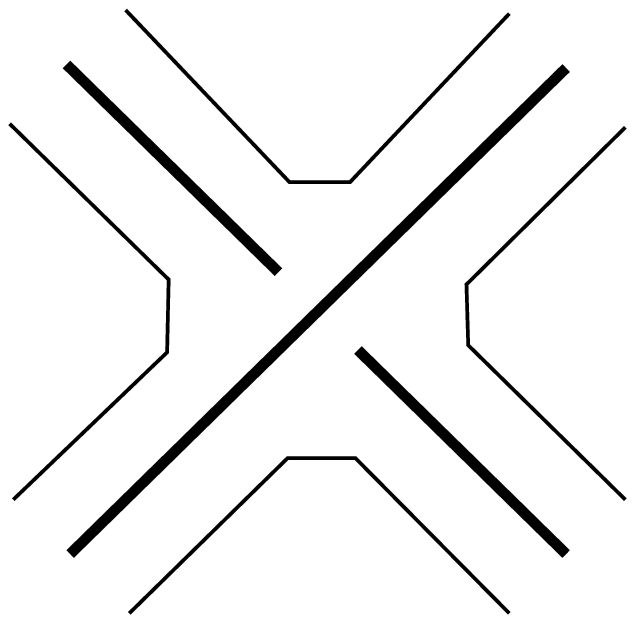}}}
\def\crsAald{\vbox to 12pt{\epsfxsize=25pt
\epsfbox{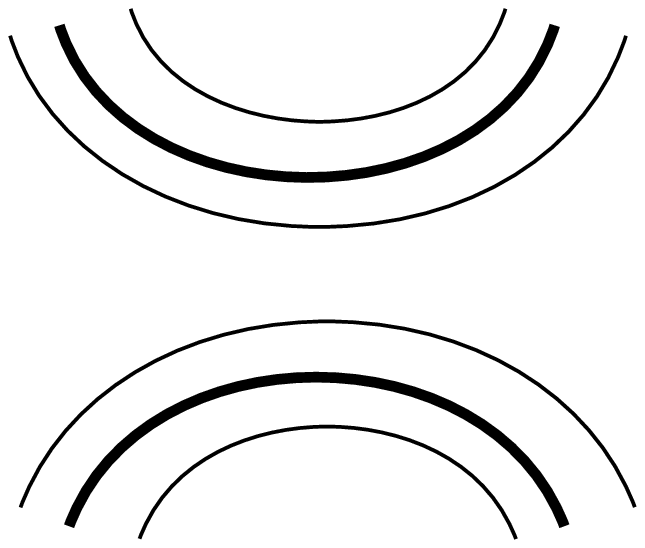}}}
\def\crsBald{\vbox to 15pt{\epsfxsize=20pt
\epsfbox{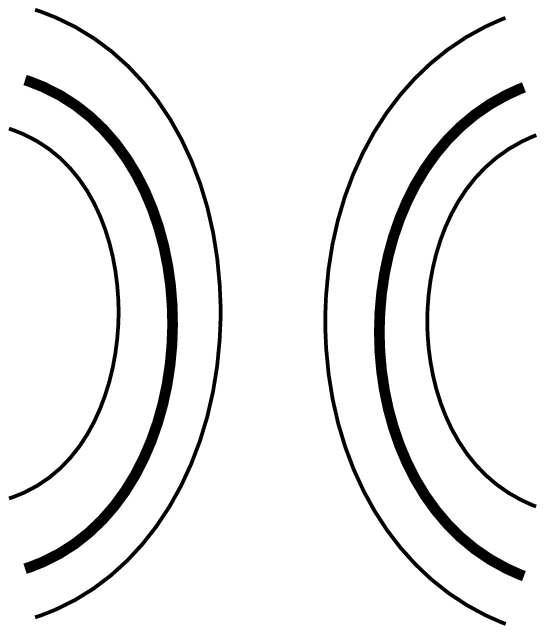}}}

\newtheorem{thm}{Theorem}
\newtheorem{cor}[thm]{Corollary}
\newtheorem{lem}[thm]{Lemma}
\newtheorem{prop}[thm]{Proposition}
\newtheorem{clm}[thm]{Claim}
\renewcommand{\theenumi}{\roman{enumi}}

\title{On the Jones polynomials of \\
checkerboard colorable virtual
knots}
\author{Naoko Kamada\\Department of Mathematics and Statistics
\\University of South Alabama}
\maketitle

\begin{abstract} 
In this paper we study 
the Jones polynomials of virtual links and abstract links. 
It is proved that a certain property of the Jones 
polynomials of classical links is valid 
for virtual links which admit checkerboard colorings.  
\end{abstract}



\section{Introduction}
In 1996, L.~H.~Kauffman introduced 
the notion of a virtual knot, which is motivated by 
study of knots in a thickened surface and abstract Gauss codes, 
cf. \cite{rkauC, rkauD}.  
M.~Goussarov, M.~Polyak, and O.~Viro \cite{rGPV} 
proved that the natural map from the category of classical knots 
to the category of virtual knots is injective; namely,  
if two classical knot diagrams are equivalent as virtual knots, 
then they are equivalent as classical knots.  
Thus, virtual knot theory is a generalization of knot theory. 
It is also found in their paper \cite{rGPV} that 
the notion of a virtual knot is helpful to study of finite type invariants.   

Kauffman defined the Jones polynomial of a virtual knot, 
which is also called the normalized bracket polynomial 
or the $f$-polynomial (cf. \cite{rkauD}).  In this paper, 
according to  \cite{rkauD}, 
we call it the $f$-polynomial instead of the Jones polynomial, 
since the definition is different from Jones' in 
\cite{rJonesA, rJonesB}.  
Finite type invariants derived from 
the $f$-polynomials are studied in \cite{rkauD}, and it is 
proved that a certain property of them 
(Corollary~14 of \cite{rkauD}) is hold in  
the category of virtual knots. 

The $f$-polynomial (Jones polynomial) of a virtual link is
quite  different from $f$-polynomials of classical links.  
For a Laurent polynomial $f$ on valuable $A$, 
we denote by ${\rm EXP}(f)$ the set of 
integers appearing as exponents of $f$. 
For example, if $f= 3A^{-2} + 6A - 7A^5$, then 
${\rm EXP}(f) = \{ -2, 1, 5\}$.   
It is well-known that for a classical link $L$ with $n$ components, 
the $f$-polynomial satisfies that 
${\rm EXP}(f) \subset 4 {\bf Z}$ if $n$ is odd, and 
${\rm EXP}(f) \subset 4 {\bf Z} +2 $ if $n$ is even.  
However, this is not true for a virtual knot/link in general. 
In this paper we introduce the notion of {\it 
checkerboard coloring} of a virtual link diagram as 
a generalization of checkerboard coloring of a 
classical link diagram.  

\begin{thm}\label{thm:checkerboard} 
Let $f$ be the $f$-polynomial of  
a virtual link $L$ with $n$ components.  
Suppose that $L$ has a virtual link diagram 
which admits a checkerboard coloring. Then 
${\rm EXP}(f) \subset 4 {\bf Z}$ if $n$ is odd, and 
${\rm EXP}(f) \subset 4 {\bf Z} +2 $ if $n$ is even.  
\end{thm}

\begin{figure}[htb] 
\begin{center}
\mbox{\epsfxsize=7cm \epsfbox{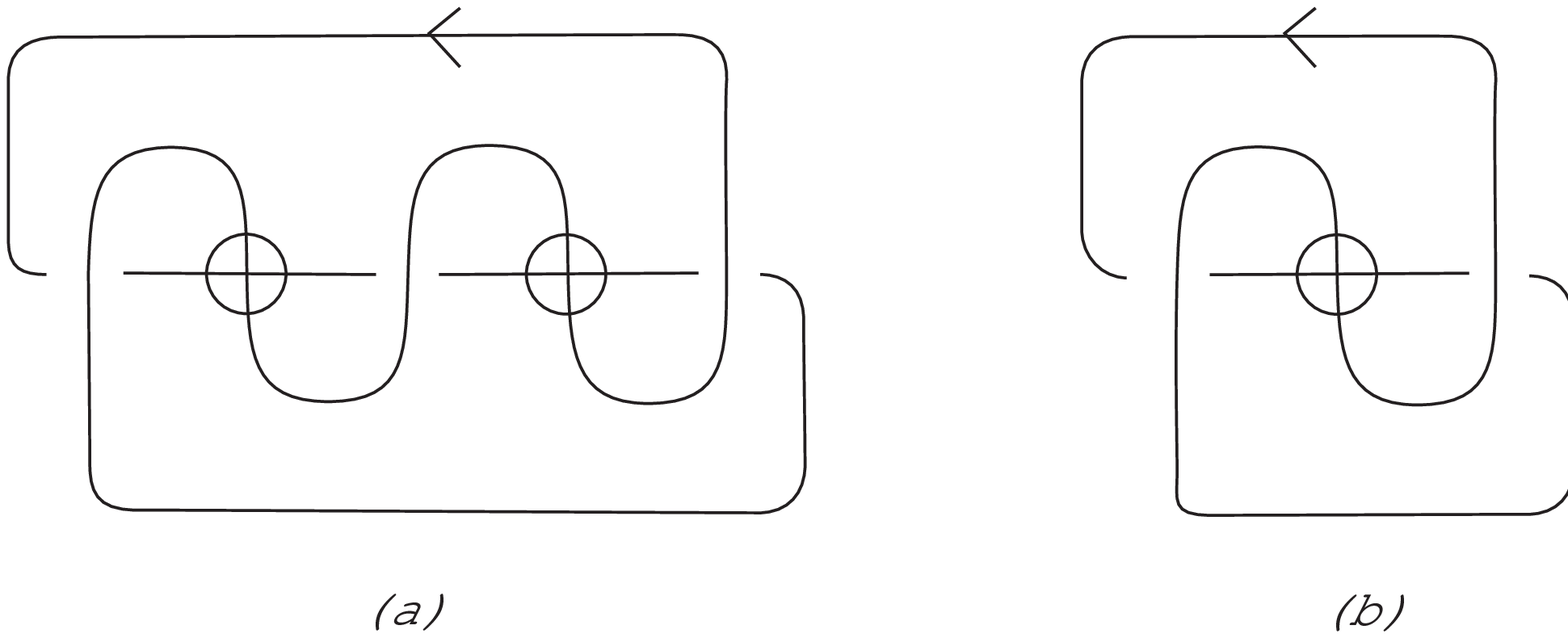}}
\end{center}
\caption[]{}\label{fig:example1}
\end{figure}

For example the virtual knot diagram illustrated 
in Figure~\ref{fig:example1} (a) 
admits a checkerboard coloring 
and the $f$-polynomial is $A^4+A^{12}-A^{16}$.  So 
${\rm EXP}(f) \subset 4 {\bf Z}$.  On the other hand, 
virtual knot diagram illustrated 
in Figure~\ref{fig:example1}
(b) does not admit a checkerboard coloring and 
the $f$-polynomial is $-A^{10}+A^6+A^4$. 
Theorem~\ref{thm:checkerboard} implies that 
this diagram is never equivalent to a 
diagram that admits a checkerboard coloring.  

%\begin{figure}[htb]
%\vspace{4cm} \hspace{1cm}
%\special{illustration example1 scaled 450}
%\caption[]{}\label{fig:example1}
%\end{figure}

If a virtual link diagram 
is alternating (the definition is 
given later), then the diagram admits a
checkerboard coloring.  Thus we have the following. 

\begin{cor}\label{cor:alternating}
Let $f$ be the $f$-polynomial of a virtual link $L$ 
with $n$ components. Suppose that $L$ has an 
alternating virtual link diagram.  
Then 
${\rm EXP}(f) \subset 4 {\bf Z}$ if $n$ is odd, and 
${\rm EXP}(f) \subset 4 {\bf Z} +2 $ if $n$ is even.  
\end{cor}

By this corollary, we see that the virtual knot 
represented by Figure~\ref{fig:example1}
(b) is not equivalent to an alternating diagram.

\section{Virtual link diagram and abstract link diagram} 

\begin{figure}[htb] 
\begin{center}
\mbox{
\epsfxsize=8cm 
\epsfbox{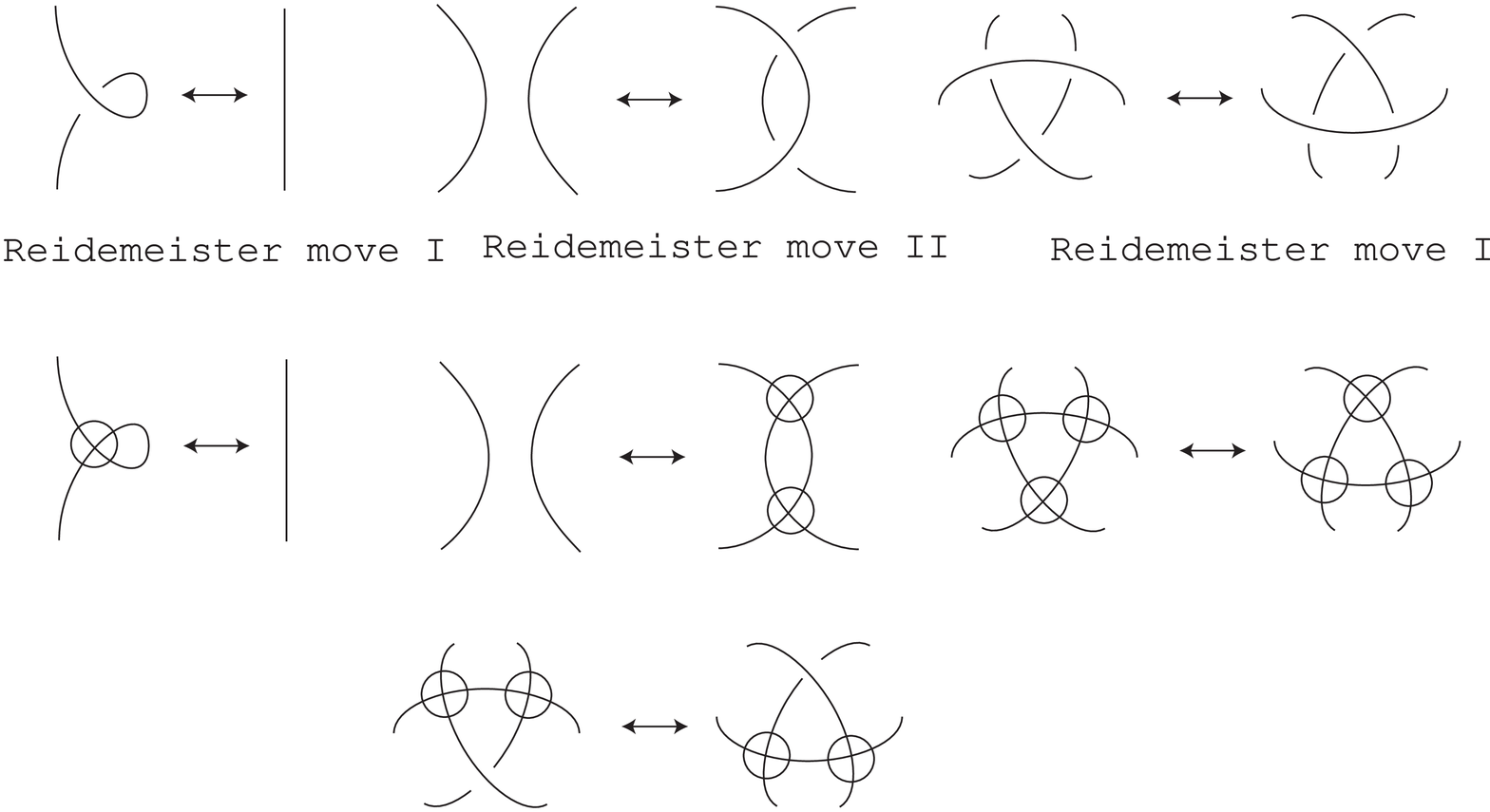}}
\end{center}
\caption[]{}\label{fig:rmovesB}
\end{figure}

A {\it virtual link diagram\/} is a closed oriented 1-manifold  
generically immersed in ${\bf R}^2$ 
such that each double point has   
information of a crossing (as in classical knot theory) 
or a virtual crossing which is indicated by a small 
circle around the double point. 
The moves of virtual link diagrams illustrated 
in Figure~\ref{fig:rmovesB} are 
called {\it generalized Reidemeister moves\/}. 
Two virtual link diagrams are 
said to be {\it equivalent\/} 
if they are related by a 
finite sequence of generalized Reidemeister moves. 
We call the equivalence class of 
a virtual link diagram {\it a virtual link\/}.

%\begin{figure}[htb]
%\vspace{5.5cm} \hspace{1cm}
%\special{illustration rmovesB scaled 300}
%\caption[]{}\label{fig:rmovesB}
%\end{figure}

\begin{figure}[htb] 
\begin{center}
\mbox{
\epsfxsize=8cm 
\epsfbox{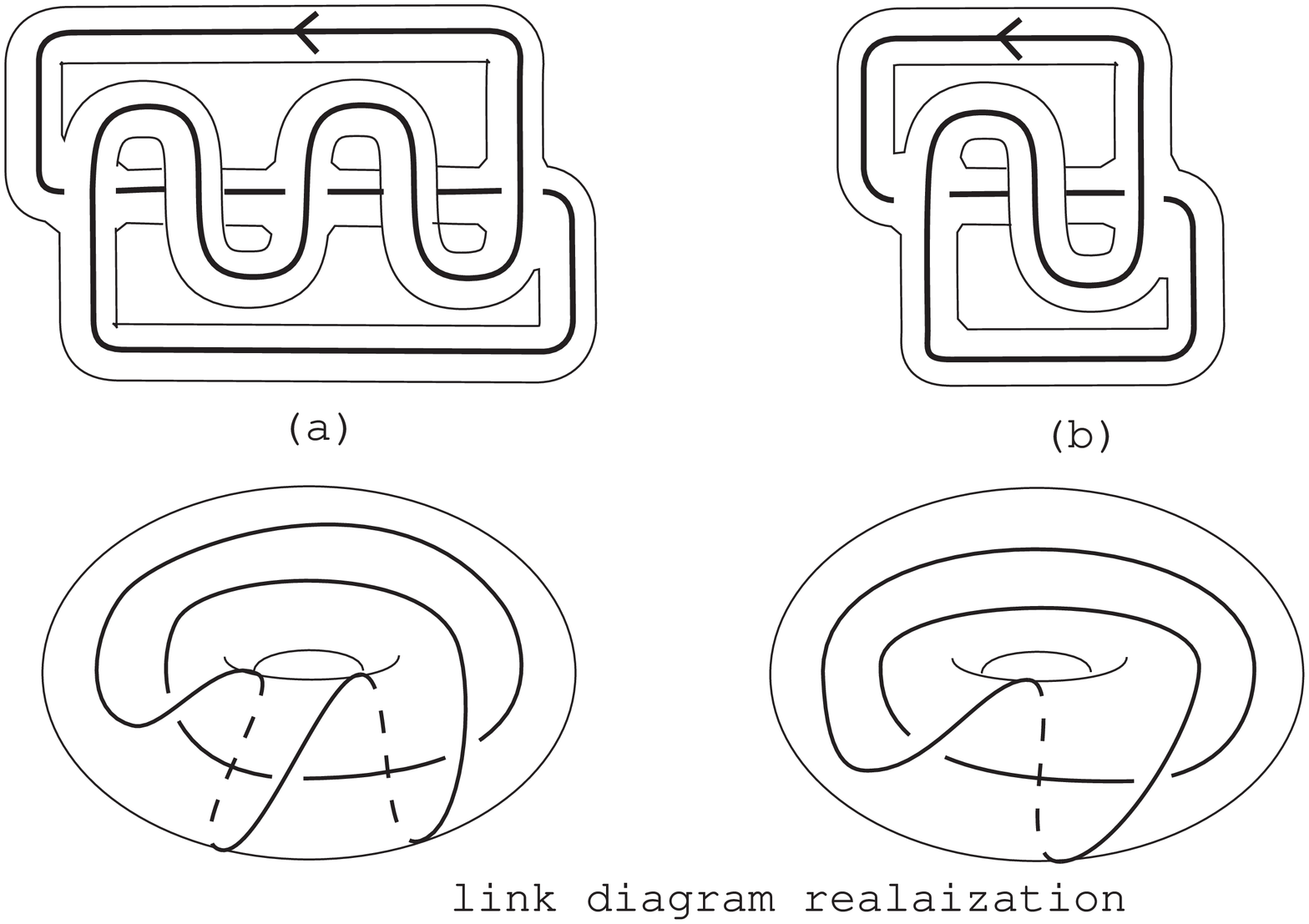}}
\end{center}
\caption[]{}\label{fig:example2}
\end{figure}

A pair $P=(\Sigma,D)$ of a compact oriented surface $\Sigma$ 
and a link diagram $D$ in $\Sigma$ is called an {\it abstract link
diagram\/} (ALD) if $|D|$ 
is a deformation retract of $\Sigma$, where
$|D|$ is a graph obtained from $D$ by replacing 
each real/virtual crossing point
with a vertex.  
For an ALD, $P=(\Sigma,D)$, if there
is an orientation preserving embedding 
$f:\Sigma\rightarrow F$ into a closed
oriented surface $F$, $f(D)$ is a link diagram in $F$. 
We call it a {\it link diagram
realization\/} of $P$ in $F$. 
In Figure~\ref{fig:example2}, we show two 
abstract link diagrams and their link diagram
realizations.  
Two ALDs  $P=(\Sigma,D)$,  $P'=(\Sigma ',D')$ are 
related by an {\it abstract Reidemeister move\/} 
(of type I, II or III)  
if there is a closed oriented surface $F$ 
and link diagram realizations of $P$ and $P'$ in $F$ 
which are related by a Reidemeister move 
(of type  I, II or III) in $F$.  
Two ALDs are {\it equivalent\/} if they are 
related by a finite sequence of abstract 
Reidemeister moves.  
We call the 
equivalence class of an ALD {\it an abstract link\/}. 

%\begin{figure}[htb]
%\vspace{7cm} \hspace{1cm}
%\special{illustration example2 scaled 350}
%\caption[]{}\label{fig:example2}
%\end{figure}

\begin{figure}[htb] 
\begin{center}
\mbox{\epsfxsize=6cm \epsfbox{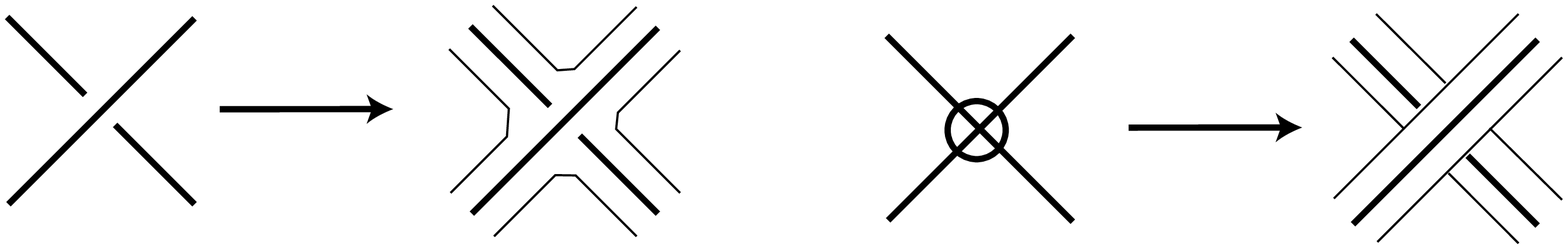}}
\end{center}
\caption[]{}\label{fig:virtualaldB}
\end{figure}

In \cite{rkk} a map 
$$\phi :\{ \text{virtual link diagrams} \} \longrightarrow
\{\text{ALDs} \}$$ 
was defined.  The idea of this map is 
illustrated in Figure~\ref{fig:virtualaldB}.  
Refer to \cite{rkk} for the definition. 
We call $\phi (D)$ an {\it ALD associated with\/} 
a virtual link diagram $D$. 
The ALDs in
Figure~\ref{fig:example2} (a) and (b)  are  ALDs 
associated with the virtual link  
diagrams in Figure~\ref{fig:example1} (a) and (b) respectively. 

%\begin{figure}[htb]
%\vspace{1.8cm} \hspace{2cm}
%\special{illustration virtualaldB scaled 300}
%\caption[]{}\label{fig:virtualaldB}
%\end{figure}

\begin{thm}(\cite{rkk})
The map $\phi$ induces a bijection
$$\Phi :\{ \text{virtual links}\} \longrightarrow 
\{ \text{abstract links} \}$$.
\end{thm}

Let $P=(\Sigma,D)$ be a pair of a compact oriented surface $\Sigma$ 
and a link diagram $D$ in $\Sigma$.  
A {\it checkerboard coloring\/} is a coloring 
of the all components of $\Sigma - |D|$ by 
two colors, say black and white, such that 
two components of  $\Sigma - |D|$ which are 
adjacent by an edge of $D$ have always distinct colors.  

We say that a virtual link diagram 
{\it admits a checkerboard coloring\/} or 
it is {\it checkerboard colorable\/} if 
the associated ALD admits a checkerboard coloring.

\section{The $f$-polynomials of abstract link diagrams}

An ALD, $P=(\Sigma, D)$, is said to be {\it unoriented\/} 
if the diagram $D$ is unoriented. 
There is a unique map 
$$<\quad>:\{\text{unoriented ALDs} \}\longrightarrow 
\Lambda = {\bf Z}[A,A^{-1}]$$ 
satisfying the following rules.

\begin{enumerate}
\renewcommand{\labelenumi}{(\theenumi)}
\item $<T>_F=1$ where $T$ is a one-component trivial ALD,
\item $<T \amalg D>=(-A^2-A^{-2})<D>$ if $D$ is not empty, 
where $\amalg$
means the disjoint union, and
\item $<\crsald>=A<\crsAald>+A^{-1}<\crsBald>$.
\end{enumerate}
Then $<\quad>$ is an invariant under 
abstract Reidemeister moves II and III. 
We call it the {\it Kauffman bracket polynomial\/} of ALD, 
cf. \cite{rkamA}.

%\begin{figure}[htb]
%\vspace{3cm} \hspace{4cm}
%\special{illustration splice scaled 300}
%\caption[]{}\label{fig:splice}
%\end{figure}

\begin{figure}[htb] 
\begin{center}
\mbox{\epsfxsize=5cm \epsfbox{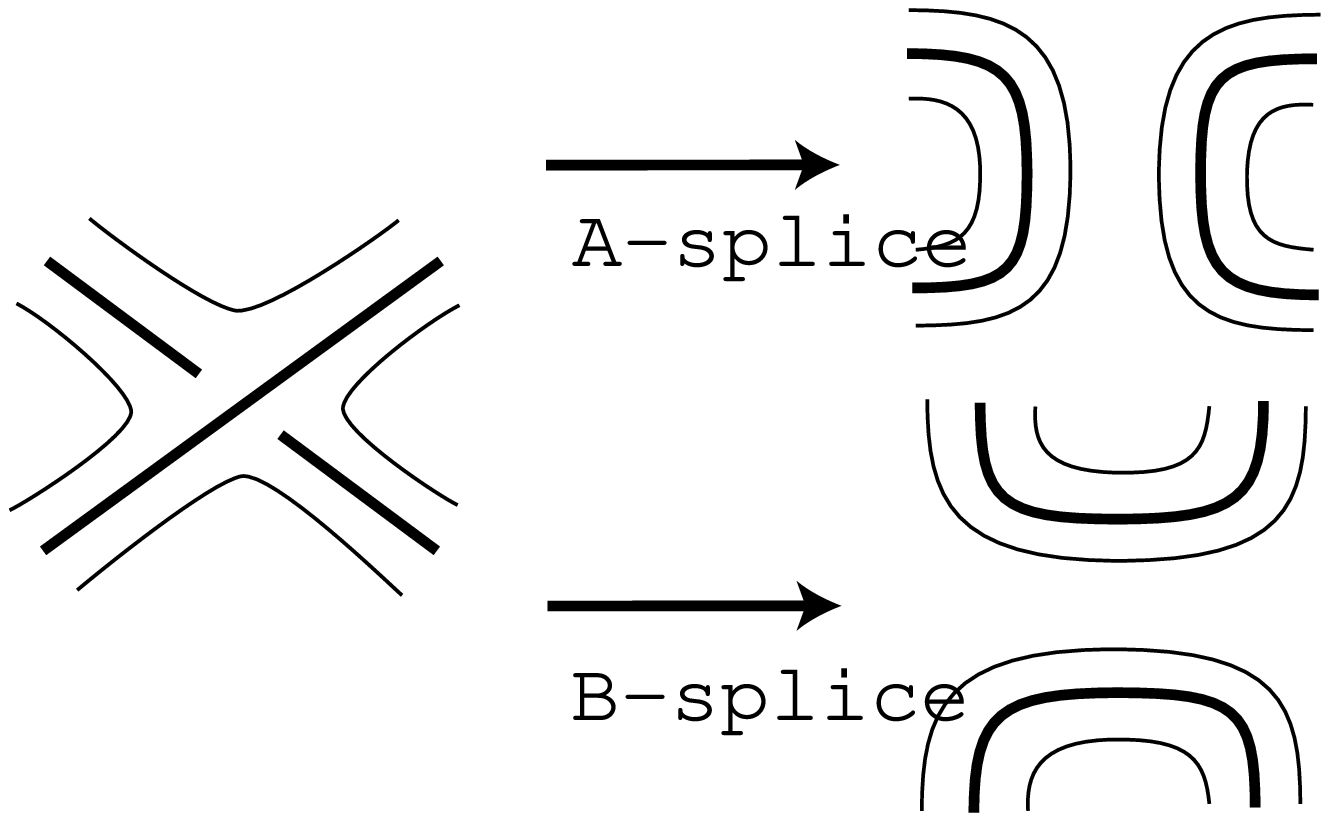}}
\end{center}
\caption[]{}\label{fig:splice}
\end{figure}

Let $P=(\Sigma, D)$ be an unoriented ALD. 
Replacing the neighborhood of a double point 
as in Figure~\ref{fig:splice}, we have
another unoriented ALD.  We call it 
an unoriented ALD obtained from $D$
by doing an {\it A-splice~} or {\it B-splice~} at 
the crossing point. 
An unoriented trivial ALD obtained from $P$ 
by doing an A-splice or B-splice at each 
crossing point is said to be a {\it state~} of $P$. 
From the definition of $<\quad>$, we see
$$<P>=\sum_{S}A^{\natural(S)}(-A^2-A^{-2})^{\sharp (S)-1},$$
where $S$ runs over all of states of $D$, 
$\natural(S)$ is the number of
A-splice minus that of B-splice used for obtaining $S$ 
and $\sharp (S)$ is the number of components of $S$. 

For an ALD, $P=(\Sigma, D)$, the writhe $\omega (P)$ is defined
by the number of positive crossings minus 
the number of negative crossings.  
Then we define the {\it normalized bracket polynomial\/} 
or the {\it $f$-polynomial\/} of $P$  
by 
$$f_P (A)=(-A^3)^{-\omega (P)}<P>.$$

By normalizing by $(-A^3)^{-\omega (P)}$, this value 
is preserved under abstract Reidemeister moves  
of type~I.  Thus this is an invariant of 
an abstract link.  
This invariant was defined in \cite{rkamA}, where 
it is called the Jones polynomial of $P$.   
It should be noted that the bijection 
$\Phi$ preserves the $f$-polynomial.

\section{Proof of Theorem~\ref{thm:checkerboard}}

Let $p$ be a crossing point of an ALD, $P=(\Sigma, D)$.  
Let $P_0 = (\Sigma_0, D_0)$ and 
$P_\infty=(\Sigma_\infty, D_\infty)$ 
be ALDs obtained from $P$ by splicing at $p$ 
orientation coherently and orientation incoherently, 
respectively.  Note that $D_\infty$ does not inherit 
an orientation from $D$.  The crossing point $p$ is either 
(i) a self-intersection 
of an immersed loop of $D$ or (ii) an intersection of two 
immersed loops.   Let $\alpha$ and $\alpha'$ be the immersed 
open arcs obtained from the loop (in case (i)) or 
from the two loops (in case (ii)) by cutting at $p$.  
Choose one of them, say $\alpha$, and we give an 
orientation to $D_\infty$ which is induced from 
that of $D$ except $\alpha$ (and hence the  
orientation is reversed on $\alpha$).  
Let $C$ be the set of crossing points of 
$D$, except $p$, such that the sign of the 
crossing point does not change in $D$ and $D_\infty$; in other word, 
at each crossing point belonging to $C$, 
both of the two intersecting arcs are contained in $D-\alpha$ 
or both of them are in $\alpha$.  
Let $C'$ be the set of crossing points of 
$D$, except $p$, such that the sign of the 
crossing point changes in $D$ and $D_\infty$; in other word, 
at each crossing point belonging to $C'$, 
one of the two intersecting arcs is contained in $D-\alpha$ 
and the other is in $\alpha$.  
Let $k$ (or $\ell$, resp.) be the number of 
positive crossings of $C$ (resp. $C'$) minus 
the number of negative crossings of $C$ (resp. $C'$).  

\begin{lem}\label{lem:jones3} 
In the above situation, 
let $f$, $f_0$ and $f_\infty$ be the 
$f$-polynomials of $P$, $P_0$ and $P_\infty$, respectively. 
Then we have 
$$f= \left\{
\begin{array}{ll}
-A^{-2} f_0 - (-A^3)^{-2 \ell} A^{-4}f_\infty, 
& \text{\rm if $p$ is a positive crossing,}\\
-A^{+2} f_0 - (-A^3)^{-2 \ell} A^{+4}f_\infty,  
& \text{\rm if $p$ is a negative crossing.} 
\end{array} \right. $$
\end{lem} 

{\it Proof.} \quad 
If $p$ is a positive crossing, then 
the writhes are $\omega(D)=k+ \ell +1$, $\omega(D_0)=k+ \ell$ and
$\omega(D_\infty)=k - \ell$.  
Since $<P>= A <P_0> + A^{-1}<P_\infty>$, we have the result. 
The case that $p$ is a negative crossing is similar. 
\qed

\bigskip

{\bf Remark.} \quad In the remark of Section~5 of
\cite{rkauD}(page~677),  an equation which is similar to 
Lemma~\ref{lem:jones3} is given.  However, it seems 
to be forgotten there to take account of the term 
$(-A^3)^{-2 \ell}$.  In consequence, the recursion formula  
of Theorem~13 of \cite{rkauD} is as follows: 
$$
v_n(G_*) = 
\sum_{k=0}^{n-1} 
\frac{2^{n-k}}{(n-k)!} 
\{ (1-(-1)^{n-k}) v_k(G_0) + 
\{ (2-3 \ell)^{n-k} - (-2 -3\ell)^{n-k} \} v_k(G_\infty) \}.$$ 
By this formula, Corollary~14 of \cite{rkauD} is still true.

\begin{cor} \label{cor:jones4} 
{\rm (cf. Theorem~13 of \cite{rkauD})} \/  
Let $f$ be the $f$-polynomial of an ALD
with $n$ components.  Then $f(1)= (-2)^{n-1}$.  
In particular, $f$-polynomials of ALDs are not zero. 
\end{cor}

{\it Proof.} \quad 
It follows from Lemma~\ref{lem:jones3} by induction 
on the number of (real) crossing points. 
\qed

\bigskip

Since $\Phi$ preserves the $f$-polynomials, 
Theorem~\ref{thm:checkerboard} is equivalent to the following 
theorem.    

\begin{thm}\label{thm:checkerboardald} 
Let $f$ be the $f$-polynomial of  
an ALD, $P=(\Sigma, D)$, with $n$ components.  
Suppose that $P$ admits a checkerboard coloring. Then 
${\rm EXP}(f) \subset 4 {\bf Z}$ if $n$ is odd, and 
${\rm EXP}(f) \subset 4 {\bf Z} +2 $ if $n$ is even.  
\end{thm}

{\it Proof.}  \quad 
For a state $S$ of $P$, we define $I(S)$ by 
$$I(S)=A^{\natural(S)}(-A^2-A^{-2})^{\sharp(S)-1}$$ 
so that the bracket polynomial of $P$ is the sum of 
$I(S)$ for all states $S$.  
Let ${\rm ind}(S)$ be a value in ${\bf Z}_4=\{0,1,2,3\}$ 
such that $I(S) \subset 4{\bf Z} + {\rm ind}(S)$.

\begin{figure}[htb] 
\begin{center}
\mbox{
\epsfxsize=6cm \epsfbox{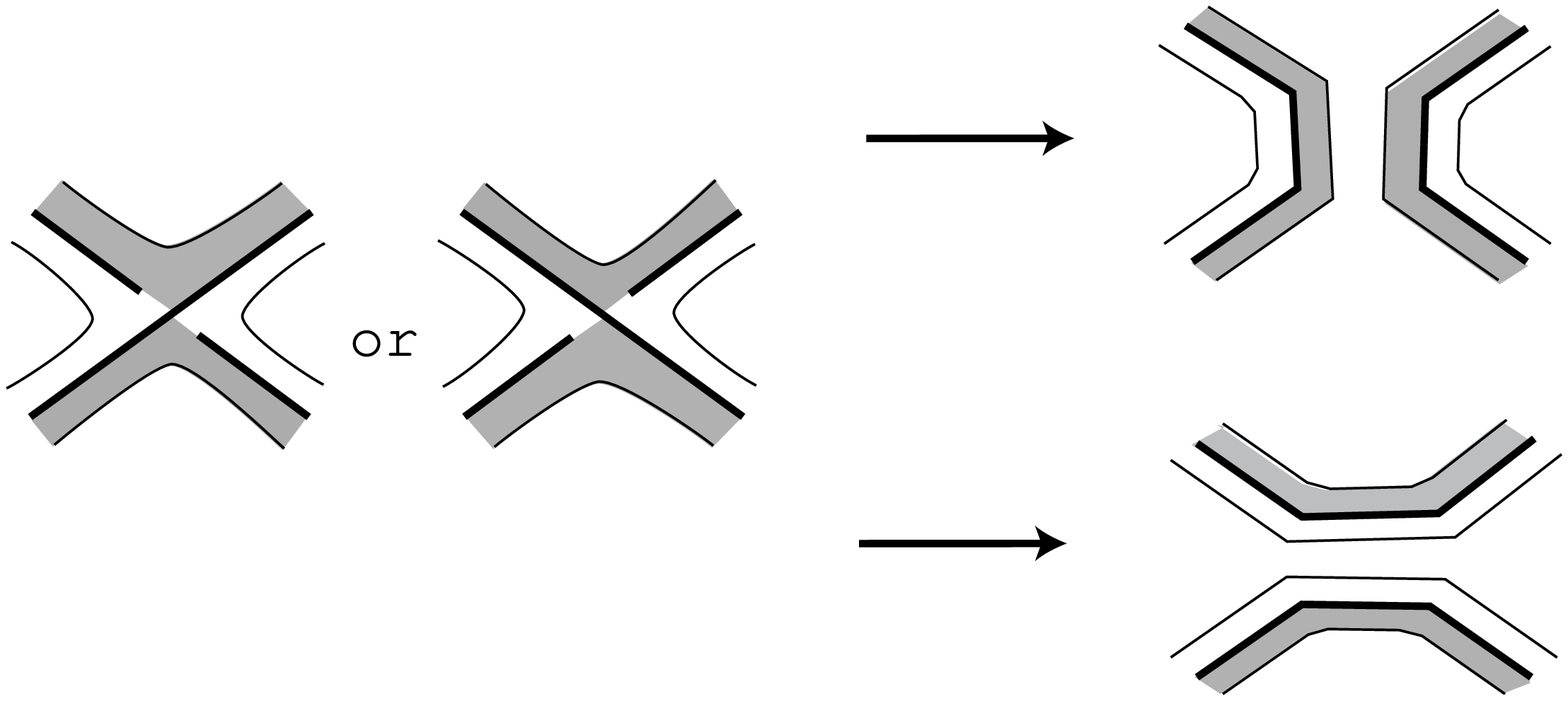}}
\end{center}
\caption[]{}\label{fig:checkersplice}
\end{figure}

\begin{figure}[htb] 
\begin{center}
\mbox{\epsfxsize=6cm \epsfbox{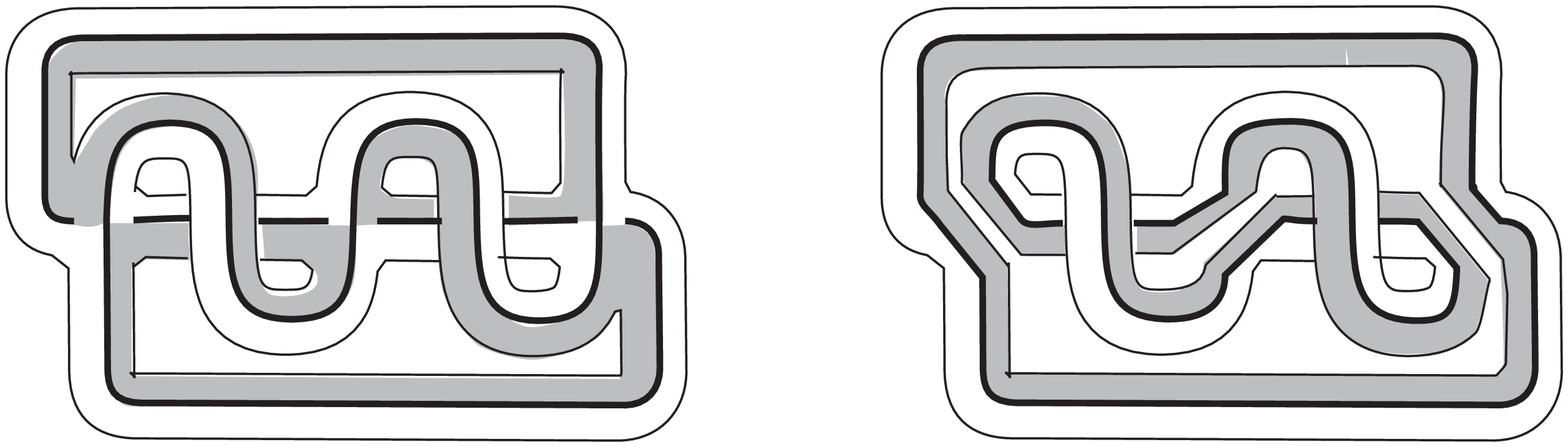}}
\end{center}
\caption[]{}\label{fig:checkstateB}
\end{figure}

\begin{figure}[htb] 
\begin{center}
\mbox{\epsfxsize=5cm \epsfbox{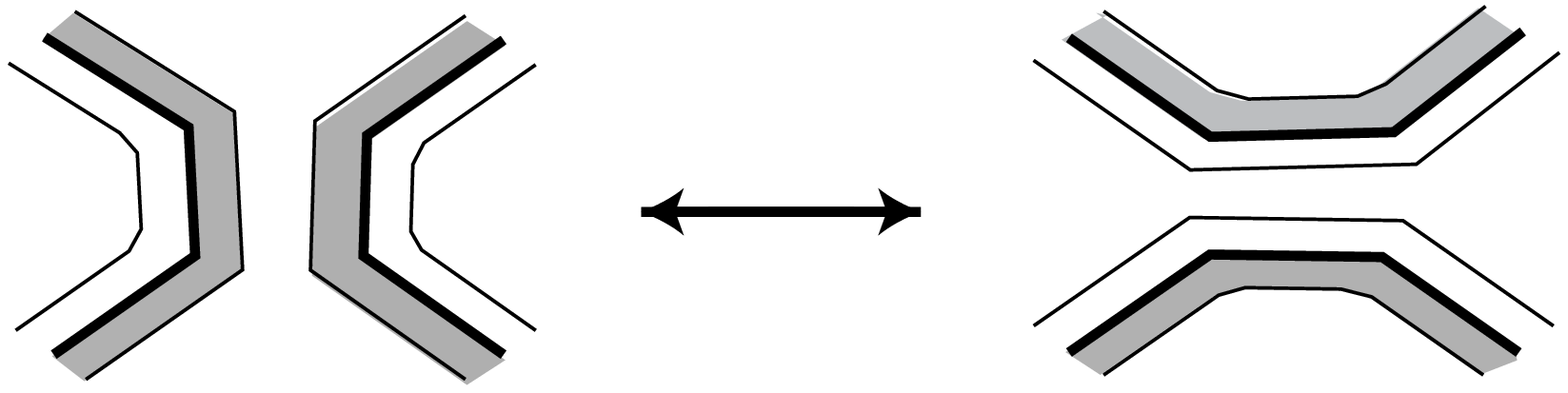}}
\end{center}
\caption[]{}\label{fig:switchspliceB}
\end{figure}

\begin{figure}[htb] 
\begin{center}
\mbox{\epsfxsize=6cm \epsfbox{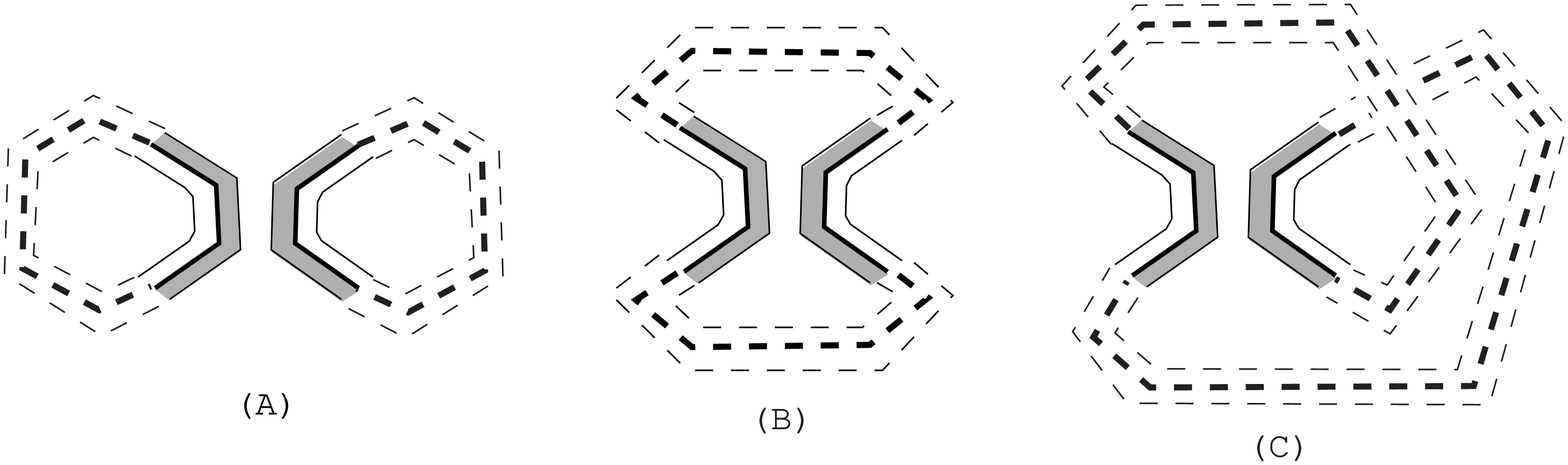}}
\end{center}
\caption[]{}\label{fig:connection}
\end{figure}

Every state of $P$ has a unique 
checkerboard coloring induced from the 
checkerboard coloring of $P$, see Figure~\ref{fig:checkersplice}. 
(Figure~\ref{fig:checkstateB} is an example of an 
ALD with a checkerboard coloring and 
a state with the induced checkerboard coloring.) 
Using this fact, we prove that ${\rm ind}(S)={\rm ind}(S')$ 
for any states $S$ and $S'$ of $P$.  
It is sufficient to prove this 
in a special case that $S$ and $S'$ are the same state 
except a crossing point, say $p$, of $D$ where 
$S$ and $S'$ are as in Figure~\ref{fig:switchspliceB}.  
For this state $S$, there are two cases (A) and (B) as in 
Figure~\ref{fig:connection}.  The case (C) does not occur, 
because a state as in (C) 
does not have a checkerboard coloring induced from 
the checkerboard coloring of $P$.  
In both cases (A) and (B),  we have  
$I(S')=A^{\natural(S)\pm 2}(-A^2-A^{-2})^{\sharp(S)-1\pm 1}$ 
and ${\rm ind}(S) = {\rm ind}(S')$.

%\begin{figure}[htb]
%\vspace{3.5cm} \hspace{1.5cm}
%\special{illustration checkersplice scaled 320}
%\caption[]{}\label{fig:checkersplice}
%\end{figure}

%\begin{figure}[htb]
%\vspace{3.0cm} \hspace{1.5cm}
%\special{illustration checkstateB scaled 320}
%\caption[]{}\label{fig:checkstateB}
%\end{figure}

%\begin{figure}[htb]
%\vspace{2cm} \hspace{3.0cm}
%\special{illustration switchspliceB scaled 320}
%\caption[]{}\label{fig:switchspliceB}
%\end{figure}

%\begin{figure}[htb]
%\vspace{2.7cm} \hspace{1cm}
%\special{illustration connection scaled 250}
%\caption[]{}\label{fig:connection}
%\end{figure}

Now we have that 
${\rm EXP}(f) \subset 4{\bold Z} + i$ where 
$i = {\rm ind}(S)$ for any state $S$ of $P$.  
We denote this number $i$ by ${\rm ind}(f)$.  
The remaining task is to prove this 
index is $0$ if $n$ is odd, and $2$ if $n$ is even.  
This is proved by induction on the 
number of (real) crossing points of $P$. 
If $P$ has no real crossing points, then 
this is obvious by the definition of the $f$-polynomial. 
If there is a crossing point, say $p$,  
apply Lemma~\ref{lem:jones3}.  
Note that $P_0$ and $P_\infty$ have checkerboard 
colorings, and 
${\rm EXP}(f_0) \subset 4{\bf Z} + {\rm ind}(f_0)$ 
and  
${\rm EXP}(f_\infty) \subset 4{\bf Z} + {\rm ind}(f_\infty)$.  
Since $f \neq 0$ and $f_0 \neq 0$ 
(Corollary~\ref{cor:jones4}), 
it follows from the equation in Lemma~\ref{lem:jones3} 
that ${\rm ind}(f) = {\rm ind}(f_0) +2 \in {\bf Z}_4$.  
The ALD $P_0$ has fewer crossing  
points than $P$ and has a checkerboard coloring.  
By induction hypothesis, 
${\rm ind}(f_0)$ is $0$ if $n'$ is odd, and $2$ if $n'$ is even, 
where $n'$ is the number of components of $P_0$.  
Since $n'= n \pm 1$, we have that 
${\rm ind}(f)$ is $0$ if $n$ is odd, and $2$ if $n$ is even. 
\qed

\section{Alternating virtual link diagrams and ALDs} 

An ALD or a virtual link diagram is {\it alternating\/} 
if we meet over and under crossing points 
alternatively when we travel along each component of 
the diagram twice.  

\begin{lem}\label{lem:altcheck}
For an ALD, $P=(\Sigma, D)$, the following conditions are
equivalent.
\begin{itemize}
\item[{\rm (i)}] By applying crossing changes, 
$P$ changes into an alternating ALD. 
\item[{\rm (ii)}] $P$ has a checkerboard coloring. 
\end{itemize}
\end{lem}  

{\it Proof of Lemma~\ref{lem:altcheck}.} \quad 
If $P$ has a checkerboard coloring, 
change each real crossing according to the coloring 
as in the most left figure of Figure~\ref{fig:checkersplice}.  
Conversely if $P$ is an alternating ALD, then 
give a checkerboard coloring near each crossing point 
as in the picture, which is extended to a 
checkerboard coloring of $P$.  
\qed

\bigskip

{\it Proof of Corollary~\ref{cor:alternating}.} \quad 
It follows from Theorem~\ref{thm:checkerboard} 
and
Lemma~\ref{lem:altcheck}. \qed 

\bigskip

\begin{figure}[htb] 
\begin{center}
\mbox{\epsfxsize=3.5cm \epsfbox{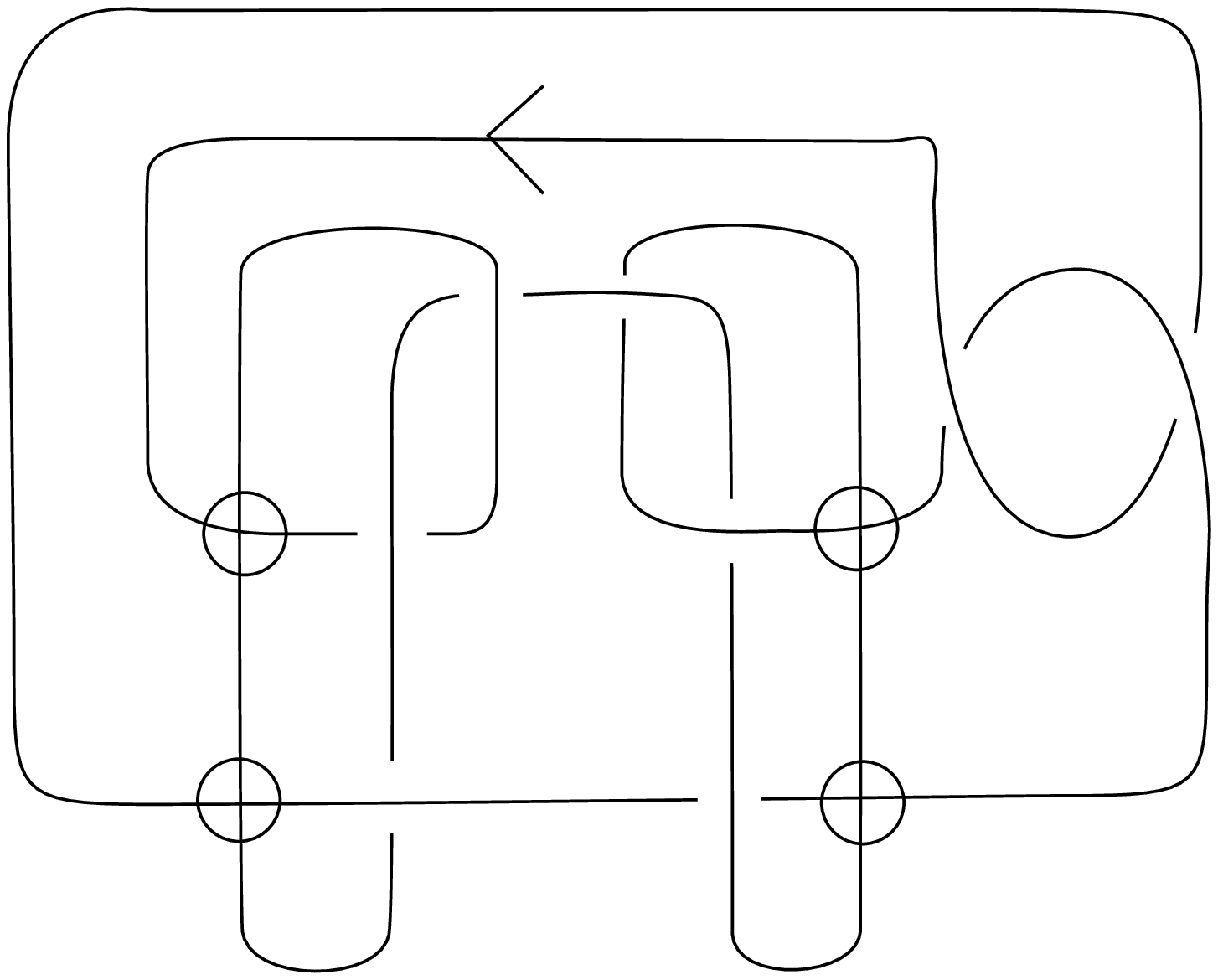}}
\end{center}
\caption[]{}\label{fig:example3}
\end{figure}

{\bf Remark.} \quad 
M.~B.~Thistlethwaite \cite{rthi} 
and K. Murasugi \cite{rmusB} 
showed that the $f$-polynomial 
(Jones polynomial) 
of a non-split alternating link is alternating, namely it
is in a form of 
$A^{\alpha}\sum c_i A^{4i}$ such that
$c_i c_j\geq 0$ for $i\equiv j\quad (\bmod\; 2)$ 
and $c_i c_j\leq 0$ for
$i\not\equiv j\quad (\bmod\; 2)$.  
This does not hold in virtual knot theory. The $f$-polynomial  
of a virtual knot in Figure~\ref{fig:example3} is
$A^{12}+3A^{16}-4A^{20}+3A^{24}-4A^{28}+4A^{32}-3A^{36}+A^{40}$.

%\begin{figure}[htb]
%\vspace{3.5cm} \hspace{3.5cm}
%\special{illustration example3 scaled 300}
%\caption[]{}\label{fig:example3}
%\end{figure}

\end{document}